\documentclass[11pt]{amsart}
\usepackage{amsmath, amsthm, amscd, amssymb, color, latexsym}
\usepackage{color,graphics,epic,epsfig}
\newtheorem{thm}{Theorem}[section]
\newtheorem*{thm*}{Theorem}
\newtheorem{cor}[thm]{Corollary}
\newtheorem*{cor*}{Corollary}
\newtheorem*{prop*}{Proposition}
\newtheorem{lem}[thm]{Lemma}
\newtheorem{prop}[thm]{Proposition}

\newtheorem*{con*}{Conjecture}
\newtheorem*{prob*}{Problem}
\newtheorem*{defn*}{Definition}

\theoremstyle{definition}
\newtheorem{defn}[thm]{Definition}
\newtheorem*{rem*}{Remark}
\theoremstyle{remark}
\newtheorem{rem}[thm]{Remark}

\newcommand{\E}{\mathcal{E}}
\newcommand{\h}{\mathcal{H}}

\newcommand{\bbR}{\mathbb{R}}

\newcommand{\bbN}{\mathbb{N}}

\begin{document}

\title{Kazhdan's Property (T) for Graphs} 

\author[Clara Brasseur]{Clara Brasseur$^{1}$}
\address{Oberlin College, 52 West Lorain Street, Oberlin OH 44074}
\email{Clara.Brasseur@oberlin.edu}
\thanks{$^{1}$Partially Supported by an N.S.F. R.E.U. Grant}

\author[Ryan E. Grady]{Ryan E. Grady$^{1}$}
\address{Clorado School of Mines, 1500 Illinois Street, Golden CO 80401}
\email{rgrady@mines.edu}

\author[Stratos Prassidis]{Stratos Prassidis$^{1,2}$}
\address{Department of Mathematics, Canisius College, 2001 Main Street, Buffalo NY 14208}
\email{prasside@canisius.edu}
\thanks{$^{2}$Partially Supported by a Canisus College Summer Research Grant }

\begin{abstract}
D. A. Kahzdan first put forth property (T) in relation to the study of discrete subgroups of Lie groups of finite co-volume.  Through a combinatorial approach, we define an analogue of property (T) for regular graphs.  We then prove the basic combinatorial and metric properties of Kazhdan groups in this context.  In particular, we use our methods to construct infinite families of expanders as in the classical case.  Finally, we consider the combinatorial analogue of the group theoretic property $(\tau)$ and prove its basic properties.
\end{abstract}

\maketitle

\section{Introduction}
Expander graphs are finite graphs that have very strong connectivity properties. They have numerous
applications in computer science and in the theory of networks.  Their existence is easily demonstrated, but explicit construction is far more difficult (\cite{Lubotzky:1994}).  Mathematically, their constructions involve methods from combinatorics, number theory, and analysis. The first construction of such a sequence was given in \cite{ma}. It was realized that most constructions of expanders involved an analytic property of groups, called property (T) which was introduced by D. A. Kazhdan (\cite{ka}).

Kahzdan first put forth property (T) in relation to the study of discrete subgroups of Lie groups of finite co-volume \cite{Lubotzky:1994}. A group satisfies Property (T) if the trivial representation is an isolated point in the space of all unitary representations of the group in the Fell topology (\cite{bhv}, \cite{Lubotzky:1994}). Equivalently, for a discrete group $\Gamma$ the Markov matrix of any $\Gamma$-invariant random walk has a {\it spectral gap} (\cite{Ollivier:2003}).

The expander property of graphs is detected by the size of the second largest eigenvalue of the 
adjacency matrix of the graph. Using this observation, a method of constructing families of expanders using Cayley graphs of infinite Kazhdan groups is described (\cite{Lubotzky:1994}, \cite{Ollivier:2003}). More specifically, if  $\{ N_i \}_{i\in I}$ is a set of finite index normal subgroups of $\Gamma$,  the set of 
Cayley graphs of 
quotient groups $\left\{\Gamma/N_i\right\}_{i\in I}$ forms a family of expanders.  The size of the second largest eigenvalue is controlled by the size of the spectral gap of $\Gamma$.

That is the starting point of the current paper. We define a {\it Kazhdan graph structure} as a pair $(X, {\Gamma})$, with $X$ a locally finite graph and $\Gamma$ a Kazhdan group acting on $X$ by graph automorphisms with a finite orbit space and finite vertex stabilizers. This definition generalizes the properties of the Cayley graph of a Kazhdan group relative to a finite set of generators while also including quasi-transitive graphs whose automorphism groups contain subgroups with property (T). This allows for new constructions of families of expanders based on graph coverings of $X$.
 
\begin{thm*}[Main Theorem]
Let $(X, {\Gamma})$ be a Kazhdan graph structure with $X$ $k$-regular. Then any family of finite $k$-regular graphs $(X_i)_{i\in I}$ covered by $X$ forms a family of expanders provided, for each 
$i\in I$,  either
\begin{itemize}
\item[(i)] $\Gamma$ is contained in $\mathrm{Cov}(X,X_i)$,  or 
\item[(ii)] $\mathrm{Cov}(X,X_i)$ is a subgroup of finite index in $\Gamma$.
\end{itemize}
Furtheromore, the expanding constant depends on the Kazhdan constant of the group $\Gamma$.
\end{thm*}

The proof of the Main Theorem follows the results on the spectral gap for Kazhdan groups given 
in  \cite{Ollivier:2003}, modified to fit our definition of Kazhdan graph structures.  Algebraically, we utilize the fact that finite extensions of Kazhdan groups are Kazhdan.

Furthermore, we establish properties of Kazhdan graph structures that are analogous to properties of Kazhdan groups. In particular,  there is no Kazhdan graph structure $(T_k, {\Gamma})$ where $T_k$ is the $k$-regular tree. Using the work of Diestel and Leader (\cite{Diestel:2001}), we construct Kazhdan graph structures $(X, {\Gamma})$  such that $X$ is not a Cayley graph.  This construction shows that quasi-transitive graphs $X$ that are constructed from Cayley graphs of Kazhdan groups $\Gamma$ provide examples of Kazhdan graph structures.

Negative kernels are important functions on topological spaces and give insight into metric properties of the space.  To this end, we consider negative kernels on $X$, where $(X,\Gamma)$ is a Kazhdan graph structure; we prove that if a kernel is $\Gamma$-invariant, then it is bounded.  Since the distance function on $X$ is $\Gamma$-invariant, the immediate implication is that the generalized roundness of $X$ is $0$ (\cite{en}, \cite{ltw}, \cite{stratos}).  Using a result from \cite{en}, we deduce that $X$ cannot be embedded isometrically into a Hilbert space.

Sometimes it is not possible to use the full strength of Kazhdan's property (T), but instead one can define a group as having the the weaker property ($\tau$) with respect to some set of finite index subgroups (\cite{lub3}).  Property ($\tau$) is useful for many of the same applications as property (T), including the construction of expanders.  Thus we extend the group theoretic definition of property ($\tau$) to graphs.  

We would like to thank Terry Bisson for several enlightening discussions and Canisius College for their hospitality.

\section{Preliminaries}

We assume the reader is familiar with basic graph theory and some geometric group theory.  (For a general reference see Woess \cite{Woess:2000} or Biggs \cite{Biggs}).  In this section we give some definitions relating to groups with Kazhdan's Property (T) and expander graphs.  Loosely speaking, expander graphs are sparse graphs with strong connectivity.  

\begin{defn}
A finite regular graph $X = X(V,E)$ with $n$ vertices and of degree $d$ is called an $(n,d,c)$-\textit{expander} if for every subset $A$ of $V$,
\begin{equation*}
\lvert \partial A \rvert \geq c\left(1-\frac{\lvert A \rvert}{n}\right)\lvert A \rvert
\end{equation*}
where $\partial A = \{ y \in V | d(y,A)=1 \}$ is the boundary of $A$ and $d$ is the distance function on $X$.  
\end{defn}

The expander property depends on the size of the second eigenvalue of the Markov matrix of the regular random walk on a graph. Also, it depends on the Cheeger constant which controls the growth of subsets inside the graph.

\begin{defn}
Let $X=X(V,E)$ be a finite graph.  Define the \textit{Cheeger constant of} $X$, denoted by $h(X)$, by:
$$h(X) = \inf_{A,B \subseteq V} \frac{\displaystyle \lvert E(A,B) \rvert}{\displaystyle \min (\lvert A \rvert , \lvert B \rvert )} $$
where the infimum runs over all the disjoint partitions $V=A \cup B$ and $E(A,B)$ is the set of edges connecting vertices in $A$ to vertices in $B$.
\end{defn}

\begin{rem}\label{rem-eigen}
It is shown in \cite{Lubotzky:1994}, Theorem 4.3.1,  that the following are equivalent 
\begin{enumerate}
 \item $\{X_i \}_{i\in I}$ is a family of $d$-regular expanders.
 \item $h(X_i) \ge {\varepsilon}_1$, for some ${\varepsilon}_1 > 0$.
\item If ${\lambda}(X_i)$ denote the largest eigenvalue of the Markov matrix on $X_i$ which is less than $1$, then ${\lambda}(X_i) \le 1 - {\varepsilon}_2$, for some ${\varepsilon}_2 > 0$.
\end{enumerate}
\end{rem}

The definitions make sense for all locally compact Lie groups but we will restrict ourselves to discrete groups.

\begin{defn}
Let $(\pi,\h)$ be a unitary representation of a  group $\Gamma$ (i.e. $\h$ is a Hilbert space and 
${\pi}: {\Gamma} \to U({\h})$ is a representation).
\begin{itemize}
\item[(i)] For a subset $Q$ of $\Gamma$ and real number $\varepsilon > 0$, a vector $\xi$ in $\h$ is $(Q,\varepsilon)$-invariant if
\begin{equation*}
\sup_{x \in Q} \| \pi (x)\xi - \xi \| < \| \xi \|.
\end{equation*}
\item[(ii)] The representation $(\pi , \h)$ has almost invariant vectors if it has $(Q,\varepsilon)$-invariant vectors for every finite subset $Q$ of $\Gamma$ and every $\varepsilon > 0$. 
\item[(iii)] The representation $(\pi,\h)$ has non-zero invariant vectors if there exists $\xi \ne 0$ in $\h$ such that $\pi ({\gamma}) \xi = \xi$ for all $\gamma \in \Gamma$.
\end{itemize}
\end{defn}

\begin{defn}
The group $\Gamma$ has \textit{Kazhdan's Property (T)} if there exists a finite subset $Q$ of $\Gamma$ and $\varepsilon >  0$ such that, whenever a unitary representation $\pi$ of $\Gamma$ has a $(Q, \varepsilon)$-invariant vector, then $\pi$ has a non-zero invariant vector (\cite{bhv}, \cite{hv}).
\end{defn}

Another important property of groups is amenability; we will see later the interaction between Kazhdan's (T) property and amenability.

\begin{defn}
A locally compact group $\Gamma$ is called \textit{amenable} if given $\epsilon > 0$ and a compact set $K \subset G$, there is a Borel set $U \subseteq G$ of positive finite (left Haar) measure $\lambda ( U )$ such that $\frac{1}{\lambda (U)} \lambda(xU \triangle U)< \epsilon$ for all $x \in K$, where $A \triangle B$ means $(A  \setminus B ) \cup ( B \setminus A )$.
\end{defn}

Amenable groups can be described in a combinatorial way: $\Gamma$ is amenable if for every $\epsilon > 0$ and every such $K$, the Cayley graph of $\Gamma$ with respect to $K$ has a finite subset, $U$, of vertices whose boundary, denoted $\partial U$, satisfies $\lvert \partial U \rvert < \epsilon \lvert U \rvert$.   (See \cite{Lubotzky:1994}).

Let $(\pi,\h )$ be a unitary linear representation of a group $\Gamma$ and a space $X$ on which a $\Gamma$-invariant reversible random walk is defined.  Following the example of \cite{Ollivier:2003} we can associate a Hilbert space to this representation. Let $\E_\pi$ be the space of $\Gamma$-equivariant functions from $X$ to $\h$, that is $f \in \E_\pi$ if for all $u \in V(X)$ and $\gamma \in \Gamma$ we have $f(\gamma u ) = \pi (\gamma) f(u)$.  We give $\E_\pi$ the inner product:
$$\langle f , g \rangle_{\E_\pi } = \sum\limits_{u \in V(X)/G } \langle f(u) , g(u) \rangle_\h ; \: f,g \in \E_\pi .
$$
Define the operator $M_\pi$ on $\E_\pi$ as follows:
$$
M_\pi (f) (u) = \frac{1}{k} \sum\limits_{v \sim u} f(v), \text{ for all } f \in \mathcal{E}_\pi .
$$

The following is the spectral characterization of property (T) (\cite{Ollivier:2003}).
\begin{thm}
The following are equivalent for a discrete group $\Gamma$:
\begin{enumerate}
\item $\Gamma$ is a Kazhdan group.
\item There is a $\sigma < 1$ such that $\text{spec}(M_{\pi}) \subset [-1, {\sigma}]{\cup}\{1\}$, for each $\Gamma$-invariant random walk.
\item There is a free random walk for which $\text{spec}(M_{\pi}) \subset [-1, {\sigma}]{\cup}\{1\}$.
\end{enumerate}
\end{thm}

\begin{rem}
Actually, in \cite{Ollivier:2003} it was shown that any $\Gamma$-invariant random walk has a two-sided spectral gap, but we will not make use of this.
\end{rem}

Using the notation from above, we can give another characterization of a group: property $(\tau)$.

\begin{defn}
Let $S$ be a finite generating set of a group $\Gamma$ and $\{N_i \}_{i \in I}$ a set of normal subgroups of finite index.  Then $\Gamma$ has {\it property} $(\tau)$ {\it with respect to} $\{N_i \}_{i \in I}$ if the Cayley graphs of each $\Gamma/N_i$ with respect to $S$ form a family of expanders i.e. the expander constant is invariant under the choice of $N_i$.   Equivalently, for each $N_i$, $M_\pi$ will have the same spectral gap for the random walk on the vertices of the Cayley graph of $\Gamma / N_i$ with respect to $S$. (This definition is equivalent to that of \cite{lub3} Proposition 2.6).
\end{defn}

Coarse equivalence or quasi-isometry between metric spaces preserves large scale invariants.

\begin{defn}
Let $(X_1 , d_1 )$ and $(X_2 , d_2 )$ be metric spaces.  We say that $X_1$ and $X_2$ are \textit{quasi-isometric} if there exist mappings
$$f: X_1 \to X_2 \: \: \text{ and }  \: \: g: X_1 \to X_2 , $$
and a constant $C>0$ such that, for all $x_1 , x_1' \in X_1$ and $x_2 , x_2' \in X_2$, we have
\begin{itemize}
\item[(i)] $\frac{1}{C} d_1 (x_1 , x_1' ) - C \le d_2 (f (x_1 ), f(x_1' )) \le C d_1 (x_1 , x_1' ) +C.$
\item[(ii)] $ \frac{1}{C} d_2 (x_2 , x_2') -C \le d_1 (g(x_2 ) , g(x_2' )) \le C d_2 (x_2 , x_2' )+C.$
\item[(iii)] $ d_1 (x_1 , g \circ f(x_1 )) \le C.$
\item[(iv)] $ d_2 (x_2 , f \circ g(x_2 )) \le C.$
\end{itemize}
\end{defn}

\begin{defn}
Given two connected graphs $X_1$, $X_2$, we say $X_2$ \textit{covers} $X_1$ if there is a map $p : X_2 \to X_1$ such that
\begin{itemize}
\item[(i)] $u \sim v$ in $X_2$ implies $p (u)  \sim p (v)$ in $X_1$.
\item[(ii)] For any $u \in X_2$, the restriction of $p$ to $N_{X_2} (u)$ is bijective onto $N_{X_1} (p (u))$.
\end{itemize}
The map $p$ is called a \textit{covering map} or simply a \textit{cover} (\cite{Woess:2000}).
\end{defn}

\begin{defn}
Let $X_1$ and $X_2$ be connected graphs such that $X_2$ covers $X_1$.  Given a covering map $p: X_2 \to X_1$, the \textit{covering transformations}, denoted $\text{Cov}(p)$, are graph automorphisms $\alpha$ such that $p \circ \alpha = p$.  In the absence of an explicit covering map we write $\text{Cov}(X_2 , X_1 )$.
\end{defn}

\begin{rem}
The general theory of covers implies that $\text{Cov}(p)$ acts freely on $V(X_2)$ and that $p$ induces a bijection between $V(X_2)/\text{Cov}(p)$ and $V(X_1)$. 
\end{rem}

\section{The Kazhdan Graph Structure}

We provide an analogue of property (T) for regular graphs; we call such graphs Kazhdan graph structures.  

\begin{defn}
A Kazhdan graph structure is a pair
 $(X, {\Gamma})$ with $X$ a regular graph and ${\Gamma} < \text{Aut}(X)$ such that
\begin{itemize}
\item[(i)] The orbit space $V(X) / \Gamma$ is finite.
\item[(ii)] For each $x \in V(X)$, the stabilizer $\Gamma_x$ is finite.
\item[(iii)] $\Gamma$ has property (T).
\end{itemize}
\end{defn}

\begin{rem}
It follows immediately from the definition that:
\begin{itemize}
\item[(i)] If $X$ is a finite graph and ${\Gamma} < \text{Aut}(X)$, then $(X, {\Gamma})$ is a Kazhdan graph structure. 
\item[(ii)] If $\Gamma$ is a Kazhdan group and $S$ is any finite symmetric  generating set of $\Gamma$, the pair $(\text{Cay}({\Gamma}, S), {\Gamma})$ is a Kazhdan graph structure.
\item[(iii)]
Let $(X, {\Gamma})$ be a Kazhdan graph structure. If $X$ is equipped with the combinatorial metric and $\Gamma$ with the word metric, then $X$ and $\Gamma$ are quasi-isometric. 
\end{itemize}
\end{rem}

Kazhdan's  property (T) is not a quasi-isometric invariant for groups (\cite{bhv}).  This observation leaves the following question open:

\vspace{12pt}\noindent
{\bf Question}. Is there a graph $X$ such that $(X, {\Gamma})$ is a Kazhdan graph structure, but $(X, {\Gamma}')$ is not (where $\Gamma$ and ${\Gamma}'$ are
two groups of automorphisms of $X$)?

\vspace{12pt}
The authors believe that such a graph exists.
\vspace{12pt}

  Our main result relates our notion of a Kazhdan graph structure and expanders.  Specifically, that for any Kazhdan graph structures $(X, {\Gamma})$, $X$ is a cover of an infinite family of expanders.  This is analogous to the result in \cite{Lubotzky:1994} for Kazhdan groups.

\begin{thm}\label{thm-main}
Let $(X, {\Gamma})$ be a Kazhdan graph structure with $X$ $k$-regular. Then the family of finite $k$-regular graphs $(X_i)_{i\in I}$ covered by $X$ forms a family of expanders provided, for each 
$i\in I$,  either
\begin{itemize}
\item[(i)] $\Gamma$ is contained in $\mathrm{Cov}(X,X_i)$,  or 
\item[(ii)] $\mathrm{Cov}(X,X_i)$ is a subgroup of finite index in $\Gamma$.
\end{itemize}
Furtheromore, the expanding constant $c$ depends on the Kazhdan constant of the group $\Gamma$.
\end{thm}

\begin{proof}
We break the proof into two cases:

\vspace{12pt}
\noindent \textsc{Case 1.} Assume the subgroup $\Gamma$ is contained in $\text{Cov}(X,X_i)$. Let
$p: X \to X_i$ be the covering map. The $\Gamma$-action on $X$ induces a $\Gamma$-action on
$X_i$ by:
$$\gamma u := p(\gamma  x),\; \text{ for } \gamma \in \Gamma, u \in V(X_i), \text{ and } x \in p^{-1} (u). $$
Since ${\Gamma} < \text{Cov}(p)$, the action is well defined. This action induces a $\Gamma$-action on the Hilbert space ${\h} = L^2(V(X_i))$ by unitary transformations. Let ${\pi}: {\Gamma} \to U({\h})$ be
the representation induced by the $\Gamma$-action.

Let $P : \h \rightarrow\h$ be the Markov operator corresponding to the regular random walk on $X_i$, defined by:
$$
P(\varphi)(u)= \frac{1}{k} \sum\limits_{v \sim u} \varphi (v), \;\; \text{for}\;\;  \varphi \in \h , \;u \in V(X_i).
$$
We define the Hilbert space $\E_\pi$ as before.  Choosing representatives $\{x_1 , x_2 , \dotsc , x_n \}$ for the orbits of $V(X_i )/ \Gamma$, the inner product on $\E_\pi$ becomes
\begin{equation*}
\langle f , g \rangle_{\mathcal{E}_\pi} = \sum\limits_{j=1}^{n} \langle f(x_j) , g(x_j) \rangle_{\h} ; \; f,g \in \mathcal{E}_\pi , 
\end{equation*}

We show that the map,
\begin{equation*}
{\iota}: \mathcal{E}_\pi \rightarrow \bigoplus_{j=1}^{n} \h,\; f \mapsto (f (x_1) , f(x_2), \dotsc f(x_n )) ,
\end{equation*} 
(where direct sum means orthogonal direct sum of Hilbert spaces) is an isomorphism. The transformation $\iota$ is obviously linear, surjective, and it preserves inner products.

\noindent
\underline{\it $\iota$ is one-to-one}: Let ${\iota}(f) = 0$. Then $f(x_j) = 0$ for $j = 1, 2, \dotsc , n$. For $u \in V(X_i)$, there is $\gamma\in \Gamma$ such that $u = {\gamma}x_j$, for some $j$. Then
$$f(u) = f({\gamma}x_j) = {\gamma}f(x_j) = 0 \;\Longrightarrow \;
f \equiv 0.$$

Recall the definition of the operator $M_\pi$ from before.  In \cite{Ollivier:2003}, it is shown that there is a number ${\sigma}_{\Gamma} < 1$, that depends only on $\Gamma$ such that $\text{spec}(M_{\pi}) \subset [-1, {\sigma}_{\Gamma}] {\cup}\{ 1\}$.

Consider the following diagram:
\begin{equation*}
\begin{CD}
\mathcal{E}_\pi    @>{\iota}>>   \bigoplus\limits_{j=1}^{n} \h \\
@V{M_\pi}VV     @VV{\bigoplus\limits_{j=1}^{n} P}V \\
\mathcal{E}_\pi @>{\iota}>> \bigoplus\limits_{j=1}^{n} \h
\end{CD}
\end{equation*}
We show that the diagram commutes, i.e.
\begin{equation*}
(\iota \circ M_\pi)(f) = \left ( \left (\bigoplus\limits_{j=1}^{n} P \right) \circ \iota \right)(f), \text{ for all } f \in \mathcal{E}_\pi.
\end{equation*} 
Indeed for $(v_1 , v_2 , \dotsc , v_n ) \in \bigoplus_{j=1}^{n} V(X_i)$ and $f \in \mathcal{E}_\pi$:
$$\begin{array}{lll}
(\iota \circ M_\pi ) (f) (v_1 , \dotsc , v_n ) &= &
 \bigl ( M_\pi (f) (x_1) , M_\pi (f)(x_2) , \dotsc , M_\pi (f)(x_n) \bigr ) (v_1 , \dotsc , v_n ) \\[2ex]
& = &\displaystyle{
\frac{1}{k} \left ( \sum\limits_{u_1 \sim v_1} f(x_1)(v_1) , \sum\limits_{u_2 \sim v_2} f(x_2)(v_2) , \dotsc , \sum\limits_{u_n \sim v_n} f(x_n)(v_n) \right )}.
\end{array}$$
Also,
$$\begin{array}{ll}
\displaystyle{
\left ( \left (\bigoplus\limits_{j=1}^{n}  P \right) \circ \iota \right)(f) (v_1 , v_2 , \dotsc , v_n )} & =  
\left ( P (f)(x_1)(v_1), P (f)(x_2)(v_2) , \dotsc , P (f)(x_n)(v_n) \right) \\[2ex]
& =  \displaystyle{
\frac{1}{k} \left ( \sum\limits_{u_1 \sim v_1} f(x_1 )(u_1 ) , \sum\limits_{u_2 \sim v_2} f(x_2)(u_2) , \dotsc , \sum\limits_{u_n \sim v_n} f(x_n)(u_n) \right )}.
\end{array}$$
Thus the spectrum of $M_{\pi}$ is equal to the spectrum of $ \oplus_{i=1}^{n} P$, which implies that the spectrum of $M_{\pi}$ equals to the spectrum of $P$, ignoring multiplicities. Therefore,
$$\text{spec}(P) \subset [-1, {\sigma}_{\Gamma}]{\cup}\{1 \}.$$

Note that, because $X_i$ is a regular graph, $1$ is an eigenvalue of  $P$.  Therefore, by Remark \ref{rem-eigen}, Part (iii)
$X_i$ is an expander where the expander constant depends only on $\Gamma$.  As our choice of $X_i$ was arbitrary, we see that $\{X_i \}_{i \in I}$ forms a family of expanders.

\vspace{12pt}
\noindent \textsc{Case 2.} \textit{Assume that $\mathrm{Cov}(X,X_i)$ is a subgroup of finite index in $\Gamma$}.  

As $\Gamma$ has property (T), $\text{Cov}(X, X_i)$ will have property (T) (cf. \cite{Ollivier:2003}).  This reduces to Case 1, where $\text{Cov}(X, X_i)$ replaces $\Gamma$.  Let $\sigma_\Gamma$ denote the Kazhdan constant with respect to $\Gamma$ and $\sigma_{\text{Cov}(X,X_i)}$ be the Kazhdan constant with respect to the group $\text{Cov}(X,X_i)$.  It follows from  \cite{Ollivier:2003} that $\sigma_{\text{Cov}(X,X_i)} \leq \sigma_\Gamma$.  Thus the expanding constant of $\{X_i\}_{i \in I}$ depends only on ${\sigma}_{\Gamma}$.
\end{proof}

\begin{rem}
With the notation as in Theorem \ref{thm-main}, if $p: X \to X_i$ is a cover with $X_i$ $k$-regular finite graph such that ${\Gamma} < \text{Cov}(p)$, then 
$$\left| V(X_i)\right| = \left| V(X)/ \text{Cov}(p)\right| < 
\left| V(X)/{\Gamma}\right|.$$
So in this case, the number of vertices of $X_i$ is bounded. Therefore, if we need to construct a sequence of expanders $\{ X_i\}_{i\in\bbN}$ so that $V(X_i) \to \infty$ we have to use covers from Case (2) in Theorem \ref{thm-main}.
\end{rem}

\section{Further Results On Kazhdan Graph Structures}

We now extend some classical results from Kazhdan groups to Kazhdan graph structures. To begin we 
characterize the Kazhdan graph structures $(X, {\Gamma})$ where $X$ is amenable. For the statement and proof corresponding to groups see \cite{Lubotzky:1994}.

\begin{prop} 
Let $(X, {\Gamma})$ be a Kazhdan graph structure with $X$ amenable.
Then $X$ is a finite graph and $\Gamma$ is a finite group.
\end{prop}

\begin{proof}
Since $X$ is amenable and the $\Gamma$-action on $X$ is quasi-transitive, $\Gamma$ is an amenable group (\cite{Woess:2000} Corollary 12.12). The only Kazhdan groups that are amenable are finite
(\cite{Lubotzky:1994}, Corollary 3.1.6), hence $\Gamma$ is finite.  As the quotient space of $\Gamma$ in $V(X)$ is finite, $V(X)$ must be finite, that is $X$ is a finite graph.
\end{proof}

It is important to understand what type of graphs appear in Kazhdan graph structures.  First, we show that there are no Kazhdan graph structures over trees.

\begin{prop}
For any natural number $k$ let $T_k$ denote the $k$-regular tree. Then, for every ${\Gamma} < \text{Aut}(T_k)$, the pair $(T_k, {\Gamma})$
is not a Kazhdan graph  structure.
\end{prop}

\begin{proof}
Suppose the contrary: that there is ${\Gamma} < \text{Aut}(T_k)$ such that $(T_k, {\Gamma})$
is a Kazhdan graph  structure. 
This means that $\Gamma$ is a uniform tree lattice.  Therefore, $\Gamma$ is finitely generated and virtually free (\cite{Bass:2001}).  As $\Gamma$ contains a free subgroup of finite index, $\Gamma$ does not have property (T); non-abelian free groups do not possess property (T) (\cite{Lubotzky:1994}), contradiction.
\end{proof}

We construct Kazhdan graph structures $(X, {\Gamma})$ such that $X$ is not a Cayley graph. 

\begin{prop}
There exist Kazhdan graph structures $(X, {\Gamma})$ so that $X$ is not a Cayley graph.
\end{prop}

\begin{proof}
Our construction is a modification of the construction of non-Cayley graphs 
in \cite{Diestel:2001} (also \cite{Watkins:1989}).  Let $\Gamma$ be a 
group with property (T) and a finite symmetric generating set $S$, e.g. 
$SL_3 (\mathbb{Z})$.  We assume that $|S| \ge 5$. Now let $G$ be the 
Cayley graph of $\Gamma$ with respect to $S$. Thus $V(G) = \Gamma$ and the
edges have the form $\{{\gamma}, {\gamma}s\}$, where ${\gamma}\in\Gamma$,
$s\in S$. Choose two different natural 
numbers $p,q$ such that $p+q = \lvert S \rvert$ with at least one of them odd. 
Without loss of generality we may assume that $q > p$.  Now let $G'$ be the 
graph obtained from $G$ by replacing each vertex by a  copy of $K_{p,q}$ (the 
complete bipartite graph with vertex classes $p$ and $q$).  
Let $K_{p,q}^{({\gamma})}$ be the copy of the bipartite graph corresponding 
to a vertex $\gamma$ of $G$.  For each edge $\{{\gamma},{\gamma}s\}$ of $G$, 
identify a vertex in the vertex class $p$ of the $K_{p,q}^{({\gamma})}$ with 
a vertex in the vertex class $q$ of the $K_{p,q}^{({\gamma}s)}$.  
We give an explicit description of the graph formed with the above 
identifications. First order the elements of $S = \{s_1, s_2, \dotsc ,
s_{p+q}\}$.
Then mark the vertices of $K_{p,q}^{({\gamma})}$ in vertex class $p$ 
by pairs $({\gamma}, s_i)$, $i = 1, \dotsc p$, 
and the elements of vertex class $q$ by 
$({\gamma}, j)$, $j = p + 1, \dotsc, p + q$.
explicitly, the vertices of $G'$ are formed from the copies of 
$K_{p,q}^{({\gamma})}$ as follows:
\begin{itemize}
\item for $i = 1, \dotsc , p$, identify the vertex $({\gamma}, s_i)$ of 
$K_{p,q}^{({\gamma})}$ to the vertex $({\gamma}s_i, s_{p+1})$ of 
$K_{p,q}^{({\gamma}s_i)}$,
\item for $j = p + 1, \dotsc , p + q$, identify the vertex $({\gamma}, s_j)$ 
of $K_{p,q}^{({\gamma})}$ to the vertex $({\gamma}s_j, s_1)$ of 
$K_{p,q}^{({\gamma}s_j)}$.
\end{itemize}

We proceed as in \cite{Diestel:2001} to show that $G'$ is not a Cayley graph.  
Assume that $G'$ is a Cayley graph.
Choose one $K_{p,q}^{({\gamma})}$ from $G'$, denote it $\widetilde{K}$.  
The graph $\widetilde{K}$ has two vertex classes, 
$\{x_1 , x_2 , \dotsc , x_p \}$ and $\{y_1 , y_2 , \dotsc , y_q \}$. 
Since $G'$ is a Cayley graph, the group of automorphisms acts transitively.
Let $\phi$ be an automorphism of $G'$ that sends an element of 
$\{y_1 , y_2 , \dotsc , y_q \}$ into $\{ y_1 , y_2 , \dotsc , y_q \}$. 
Then $\phi$ must map $\{ x_1 , x_2 , \dotsc , x_p \}$  to itself as 
$\{ x_1 , x_2 , \dotsc , x_p \}$ is the only set of vertices with $q$ 
common neighbors and having a common neighbor in 
$\{y_1 , y_2 , \dotsc , y_q \}$. Thus $\phi$ maps $\widetilde{K}$ to itself. 
Since at least one of $p$ and $q$ is odd, $\phi$ must fix a vertex in 
$\widetilde{K}$. So $G'$ does not admit a free action that is transitive 
and thus, it cannot be a Cayley graph.

We now describe an action of $\Gamma$ on $G'$. We have:
$$V(G') = \bigcup_{\gamma\in\Gamma}\bigcup_{i = 1}^{p+q}({\gamma}, s_i).$$
An element ${\gamma}'\in\Gamma$ acts by left translations: 
${\gamma}'({\gamma}, s_i) = ({\gamma}'{\gamma}, s_i)$. The action is a free 
action by graph automorphisms and the quotient $V(G')/{\Gamma} = V(K_{p,q})$. 
Since $\Gamma$ has property (T), $(G, {\Gamma})$ is a Kazhdan graph structure.
\end{proof}

\section{Kazhdan Graph Structures and Negative Kernels}

The existence of negative kernels on a topological space is an important characterization of the space; it is also of note as to whether these kernels are bounded.  It is known that any invariant negative kernel on a Kazhdan group is bounded (\cite{bhv}).  Motivated by this result we show that if $(X, \Gamma)$ is a Kazhdan graph structure then any $\Gamma$-invariant negative kernel on $X$ is bounded.  Throughout this section we will assume that $X$ is a metric space equipped with the edge-path metric.

\begin{defn}
A negative kernel on a topological space $Y$ is a continuous function $h : Y \times Y \rightarrow \mathbb{R}$ with the following properties:
\begin{itemize}
\item[(i)] $h(y,y) = 0$ for all $y$ in $Y$.
\item[(ii)] $h(y,z) = h(z,y)$ for all $y,z$ in $Y$.
\item[(iii)] For any $n$ in $\mathbb{N}$, any elements $y_1 , \dotsc , y_n$ in $Y$, and any real numbers $c_1 , \dotsc , c_n$ with$\sum_{i=1}^{n} c_i =0$, the following inequality holds:
\begin{equation*}
\sum\limits_{i=1}^{n} \sum\limits_{j=1}^{n} c_i c_j h (y_i , y_j ) \leq 0 .
\end{equation*}
\end{itemize}
\end{defn}

\begin{lem}\label{lem-negative}
Let $h: X{\times}X \to \bbR$ be a negative kernel on $X$. For each
$x\in V(X)$, define
$$h_x: {\Gamma}{\times}{\Gamma} \to {\bbR}, \;\;
h_x({\gamma}_1, {\gamma}_2) = h({\gamma}_1x, {\gamma}_2x).$$ 
Then $h_x$ is a negative kernel on $\Gamma$ for each $x\in V(X)$, when 
$\gamma$ is equipped with the word metric.
\end{lem}

\begin{proof}
We show that $h_x$ is a kernel conditionally of negative type on $\Gamma$.  Indeed, if $\gamma_1 , \gamma_2 \in \Gamma$:
\begin{itemize}
\item[(i)] 
$ h_x (\gamma_1 , \gamma_1 ) 
=  h( \gamma_1 x , \gamma_1 x ) 
=  0$.
\item[(ii)]
$h_x (\gamma_1 , \gamma_2) 
=  h(\gamma_1 x , \gamma_2 x ) 
= h(\gamma_2 x , \gamma_1 x )
= h_x (\gamma_2 , \gamma_1 )$.
\item[(iii)]
Assume $c_1 , c_2 , \dotsc , c_n$ are real numbers such that $\sum_{i=1}^{n} c_i = 0$ and $\gamma_1 , \gamma_2 , \dotsc , \gamma_n$ are elements of $\Gamma$. Then
$$
\sum\limits_{i=1}^n \sum\limits_{j=1}^n c_i c_j \left [ h_x (\gamma_i , \gamma_j ) \right ] =  \sum\limits_{i=1}^n \sum\limits_{j=1}^n c_i c_j \left [ h(\gamma_i x , \gamma_j x) \right ] 
\leq  0.$$
\end{itemize}
\end{proof}

Fix a finite set  $\{x_1, x_2, \dotsc , x_n\}$ of orbit representatives of the $\Gamma$ action on $V(X)$ and
using the above notation set $h_{x_i} = h_i$.

\begin{prop}\label{prop-bounded}
If $(X,\Gamma )$ is a Kazhdan graph structure such that $h_i$ is a bounded negative kernel on $\Gamma$ for  $i = 1, 2, \cdots n$.
Then any $h$ is bounded.
\end{prop}

\begin{proof}
Let $K$ be the maximum of the bounds of $h_i$. Also set
$$L = \max\{h(x_i, x_j):\; i, j = 1, 2, \dotsc , n\}.$$
We show that $h$ must aso be bounded by estimating a generic $h(\gamma_i x_i, \gamma_j x_j)$.

We begin by estimating $h(\gamma_j x_j, x_i)$. Let $c_1=c_2=1$ and $c_3=-2$. The fact that $h$ is a negative kernel implies that:
$$c_1 c_2 h({\gamma}_jx_j, x_i) + c_1 c_3 h({\gamma}_jx_j, x_j) + c_2 c_3 h(x_i, x_j) \le 0.$$
Therefore,
$$h(\gamma_j x_j,x_i) -2 h(\gamma_j x_j,x_j) -2 h(x_i,x_j) \le 0 \;\Longleftrightarrow\;
h(\gamma_j x_j,x_i) \le 2(h(\gamma_j x_j,x_j) + h(x_i,x_j)).$$
Thus $h(\gamma_j x_j,x_i) \le 2(K + L)$.

We apply the negative kernel property for $c_1=c_2=1$, $c_3=-2$, and $c_4=0$.
$$c_1 c_2 h({\gamma}_ix_i, {\gamma}_jx_j) + c_1 c_3 h({\gamma}_ix_i, x_i) +c_1c_4 h({\gamma}_ix_i, x_j)+ c_2 c_3 h({\gamma}_jx_j,x_i) +c_2 c_4 h({\gamma}_jx_j,x_j)+c_3 c_4 h(x_i,x_j) \le 0.$$
Therefore,
\begin{equation*}
h(\gamma_i x_i,\gamma_j x_j) -2h(\gamma_i x_i,x_i) -2 h(\gamma_j x_j,x_i) \le 0 \;\Longleftrightarrow\;
h(\gamma_i x_i,\gamma_j x_j)  \le 2(h(\gamma_i x_i,x_i) + h(\gamma_j x_j,x_i)).
\end{equation*}
Using our previous calculations and the fact that $h(\gamma_i x_i,x_i) \leq K$ we find:
\begin{equation*}
h(\gamma_i x_i,\gamma_j x_j) \le 2(K+2(K+L))=6K+4L.
\end{equation*}
Thus $h$ is bounded.
\end{proof}

A negative kernel $k$ on $X$ is called {\it $\Gamma$-invariant} if 
$$h({\gamma}x, {\gamma}x') = h(x, x'), \;\;\text{for all}\; x, x' \in X, \; {\gamma} \in \Gamma .$$

\begin{thm}\label{thm-invariant}
Let $h$ be a $\Gamma$-invariant negative kernel on $X$. Then $h$ is bounded.
\end{thm}

\begin{proof}
For each $x$ in $X$, define
$${\phi}_x : {\Gamma} \to \bbR , \;\; {\phi}_x({\gamma}) = h_x({\gamma}, 1).$$
Then ${\phi}_x({\gamma}_1^{-1}{\gamma}) = h_x({\gamma}, {\gamma}_1)$ is a negative kernel on $\Gamma$ (Lemma \ref{lem-negative}). Since $\Gamma$ is Kazhdan, ${\phi}_x$ is bounded (\cite{de}). 
The result follows from Proposition \ref{prop-bounded}.
\end{proof}

The generalized roundness of a metric space $(X,d)$ is the supremum of all $q$ such that for every $n \geq 2$ any collection of $2n$ points $\{ a_i \}^{n}_{i=1}$, $\{b_i \}_{i=1}^{n}$, we have that:
\begin{equation*}
\sum\limits_{1 \leq i < j \leq n} \left ( d(a_i , a_j )^q + d ( b_i , b_j )^q \right ) \leq \sum\limits_{1\leq i,j \leq n} d (a_i , b_j )^q .
\end{equation*}
(\cite{en}, \cite{ltw}, \cite{stratos}).

\begin{cor}
Let $(X, \Gamma )$ be a Kazhdan graph structure with $X$ infinite.  Then the generalized roundness of $X$ is $0$.
\end{cor}

\begin{proof}
Since $\Gamma$ acts by graph automorphisms, it acts by isometries on $X$. Thus the $p^{th}$ power of the distance function $d^p$ is a $\Gamma$-invariant function. If the generalized roundness of $X$ were equal to $p > 0$, then $d^p$ would be a negative kernel on $X$ (\cite{ltw}). But that would imply that $d^p$ is bounded (\ref{thm-invariant}), and thus $X$ must be finite, a contradiction.
\end{proof}

In \cite{en}, it was shown that Hilbert spaces have generalized roundness equal to $2$. Using this result, we get the following.

\begin{cor}
Let $(X, \Gamma )$ be a Kazhdan graph structure with $X$ infinite. Then there is no isometric embedding of $X$ into a Hilbert space.
\end{cor}

\section{On Property $(\tau)$}

Again we extend a group theoretic definition to graphs; in this case it is property $(\tau)$.

\begin{defn}
Let $X$ be a graph and $\{X_i \}_{i \in I}$ a sequence of finite graphs.  We say $X$ has {\it property $(\tau)$ with respect to} $\{X_i \}_{i \in I}$ if
\begin{itemize}
\item[(i)] For each $i \in I$, there exists a covering map $p_i : X \to X_i$.
\item[(ii)] $\{ X_i \}_{i \in I}$ is a family of expanders.
\end{itemize}
\end{defn}

Clearly a group with property (T) has property $(\tau)$, we show the same implication exists between our graph analogues.

\begin{prop}
Let $(X, \Gamma)$ be a Kazhdan graph structure and $\{ X_i \}_{i=1}$ a sequence of finite $k$-regular graphs covered by $X$, such that $\mathrm{Cov}(X,X_i)$ is a subgroup of finite index in $\Gamma$.  Then $X$ has property $(\tau)$ with respect to $\{X_i \}_{i \in I}$.
\end{prop}

\begin{proof}
The result follows immediately from Theorem \ref{thm-main}.
\end{proof}

The next result is the analogue of Proposition 1.32 in \cite{lub3}. We replace the residually finite property by a property on covering transformations.

\begin{prop}
Let $X$ be a graph with property $(\tau)$ with respect to $\{X_i \}_{i \in I}$ such that
$$\bigcap_{i \in I} \mathrm{Cov}(p_i ) = \{1\}.$$
If $X$ is amenable then $X$ is finite.
\end{prop}

\begin{proof}
Suppose on the contrary that $X$ is infinite. We consider a Folner sequence $\{F_n\}$ of finite subsets of $X$ such that (\cite{lub3}, Proposition 1.32)
$$\lim_{n \to \infty} \frac{\displaystyle \lvert \partial F_n \rvert}{\displaystyle \lvert F_n \rvert} =0 .$$
For each $n$, there is $j_n\in I$ such that:
\begin{itemize}
\item[(i)] $p_{j_n}|: F_n \to V(X_{j_n})$ is an injection.
\item[(ii)] $|F_n| \le |V(X_{j_n})|/2$.
\end{itemize}
 For the Cheeger constant, 
$$ h(X_{j_n} ) \le \frac{\displaystyle \lvert \partial F_n \rvert}{\lvert F_n \rvert }.$$
Thus
$$\lim_{n\to \infty} h(X_{j_n} ) = 0$$
contradicting Remark \ref{rem-eigen}, Part (ii).
\end{proof}

\end{document}